\documentclass[12pt,sn-mathphys,Numbered]{sn-jnl}

\usepackage{graphicx}%
\usepackage{multirow}%
\usepackage{amsmath,amssymb,amsfonts}%
\usepackage{amsthm}%
\usepackage{mathrsfs}%
\usepackage[title]{appendix}%
\usepackage{xcolor}%
\usepackage{textcomp}%
\usepackage{manyfoot}%
\usepackage{booktabs}%
\usepackage{algorithm}%
\usepackage{algorithmicx}%
\usepackage{algpseudocode}%
\usepackage{listings}%

\theoremstyle{thmstyleone}%
\newtheorem{thm}{Theorem}

\newtheorem{lem}[thm]{Lemma}

\theoremstyle{thmstyletwo}%
\newtheorem{rem}{Remark}%

\theoremstyle{thmstylethree}%
\newtheorem{defn}{Definition}%

\begin{document}

\title[Gradient estimates of the Finslerian Allen-Cahn equation]{Gradient estimates of the Finslerian Allen-Cahn equation}


\author*[1]{\fnm{Bin} \sur{Shen}}\email{shenbin@seu.edu.cn}



\affil*[1]{\orgdiv{School of Mathematics}, \orgname{Southeast University}, \orgaddress{\street{Dongnandaxue Road 2}, \city{Nanjing}, \postcode{211189}, \state{Jiangsu}, \country{China}}}




\abstract{In this manuscript, we study bounded positive solutions to the Finslerian Allen-Cahn equation. The Allen-Cahn equation is widely applied and connected to many mathematical branches. We find the Finslerian Allen-Cahn equation is also an Euler-Lagrange equation to a Liapunov entropy functional. Adopting the $\epsilon$-trick, we prove the global gradient estimates of its positive solutions on compact Finsler metric measure spaces adopting the lower bounds of the weighted Ricci curvature. Moreover, as application of a new comparison theorem developed by the author, we also get a local gradient estimates on noncompact forward complete Finsler metric measure spaces with locally finite misalignment, combining with the bounds of some non-Riemannian curvatures and the lower bound mixed weighted Ricci curvature. At last, we give a Liouville type theorem of such solutions.}




\maketitle

\section{Introduction}

The famous Allen-Cahn equation is an equation of motion proposed by S. Allen and J. Cahn in 1979 to describe the reversed-phase boundary in ferroalloys \cite{AC1979}, and has been widely studied as a phenomenological model of phase separation in binary alloys. At present, it has been widely used in many complicated moving interface problems in fluid dynamics, materials science, crystal growth, image inpainting and segmentation, minimal surfaces, dynamic systems, and other theoretical problems (cf. \cite{SZ2023}\cite{LTZ2020}\cite{CG2023}\cite{CM2020} etc.).

The elliptic Allen-Cahn equation and its parabolic counterpart have been extensibly studied in the last decades for their relations with the theory of minimal hypersurfaces on $\mathbb{R}^n$, with an important source of motivation being De Giorgi's conjecture \cite{Degio1979}.\\
\textbf{Conjecture. }{\it  
	Let $u\in C^2(R^n)$ be a solution to the Allen-Cahn equation
	\begin{eqnarray}\label{1.1}
		-\Delta u=u-u^3\quad \mbox{in} \quad \mathbb{R}^n
	\end{eqnarray}
	satisfying $\partial_{x^n}u > 0$. If $n \leq 8$, all level sets $\{u=\lambda\}$ of $u$ must be hyperplanes.
}

For dimension $n=2$, it was proved by N. Ghoussoub and
C. Gui \cite{GG1998}. For $n=3$, it was proved by L. Ambrosio and X. Cabr\'e \cite{AC2000}. A celebrated result by O. Savin
\cite{Savin2009} established its validity for $4\leq n\leq8$ under an extra assumption that
\begin{eqnarray*}
	\lim_{x_n\rightarrow\pm\infty}u(x',x^n)=\pm1.
\end{eqnarray*}

On the other hand, del Pino, Kowalczyk, and Wei constructed a counterexample in dimensions $n\geq 9$ \cite{PKW2011}.

The RHS of the Allen-Cahn equation is a first-order derivative of a double-well potential $W$. An intrinsic property of the Allen-Cahn equation is that it is the Euler-Lagrange equation for the following Liapunov entropy functional in the $L^2$ function space. That is,
\begin{eqnarray}
	J(u)=\int_U\left(\frac12|\nabla u|^2+W(u)\right)d\mu, \quad \mbox{with}\quad W(u)=\frac14(u^2-1)^2.
\end{eqnarray}

Solutions to the Allen-Cahn equation have the intricate connection to the minimal surface theory. 
A well known question of Yau asks that ``do all 3-manifolds contain infinite immersed minimal surfaces?" \cite{Yau1982}. Almgren and Pitts \cite{Pit1981} have developed a far-reaching theory of existence and regularity (cf. \cite{SS1981}) of min-max (unstable) minimal hypersurfaces. In particular, they proved that in a closed Riemannian manifold $(M^{n}, g)$, with $n\geq3$, there exists a closed minimal hypersurface which is smooth and embedded outside a singular set of Hausdorff dimension at most $n-8$. Such a minimal hypersurface is called a minimal hypersurface with optimal regularity.  Guaraco introduced a new approach for the min-max construction of minimal hypersurfaces which was further developed by Gaspar and Guaraco \cite{GG2018}, based on the study of the limiting behaviour of solutions to the Allen-Cahn equation. They showed that the Allen-Cahn widths are greater than or equal to the Almgren-Pitts widths. Later, A. Dey proved that the Almgren-Pitts widths and the Allen-Cahn widths coincide \cite{Dey2022}.

The gradient estimate method was successfully employed by S.T. Yau to study harmonic functions on complete Riemann manifolds \cite{Yau1975}. S. Hou proved the local gradient estimates for bounded positive solutions to \eqref{1.1} on a complete noncompact Riemann manifold \cite{Hou2018}. Later, L. Zhang established the Hessian estimates of positive solutions to \eqref{1.1} \cite{Zhang2021}. N. Dung and N. Khanh derived the gradient estimates of positive solutions to Allen-Cahn equation on weighted Riemann manifolds \cite{NN2022}.

As a generalization of such method, the gradient estimation on Finsler metric measure spaces was first given by S. Ohta to the heat equation on compact Finsler manifolds \cite{Ohta2009}. C. Xia studied the harmonic function on both compact and forward complete Finsler manifolds by adopting the Moser's iteration \cite{Xia2014}. Later, Q. Xia provided the gradient estimates of positive solutions to the heat equation on forward complete Finsler manifolds \cite{XiaQ2020} with the same technique of \cite{Xia2014}. The author proposed a generic approach to derive the global and local gradient estimates on Finsler metric measure spaces for the Finslerian Schr\"odinger equation \cite{Shen2024}.

In this manuscript, we consider the bounded positive solutions to the Finslerian Allen-Cahn equation
\begin{eqnarray}\label{equ-FAC}
	\Delta^{\nabla u}u+(1-u^2)u=0,
\end{eqnarray}
on a Finsler metric measure space $(M,F,\mu)$, where $\Delta^{\nabla u}$ is given in \eqref{def-WLap} below. We establish the gradient estimate of bounded positive solution of the Finslerian Allen-Cahn equation on both compact and forward complete Finslerian metric measure spaces following the method in \cite{Shen2024}. Concretely, we prove that

\begin{thm}\label{thm-cpt-reduced}
	Let $(M^n,F,\mu)$ be a compact Finsler $CD(-K,N)$ space, with dimension $n\geq 2$. Assume the reversibility of $M$ has upper bound $\rho_0$. Suppose $u$ is a bounded positive solution to \eqref{equ-FAC} on $M$ with $u\leq C$ where $C$ is a positive constant. 
	
	(1) If $C\leq 1$, then we have 
	$$\frac{F^2(\nabla u)}{u^2}+\frac{2}{3}(1-u^2)\leq \frac{4(N+n)}{3}\rho_0^2K.$$
	
	(2) If $C\geq 1$, then we have 
	$$\frac{F^2(\nabla u)}{u^2}+2(1-u^2)\leq 2(N+n)\rho_0^2K+27\sqrt2(N+n)^{\frac{3}{2}}C^2.$$
\end{thm}

In noncompact case, we denote the forword distance function from a fixed point $p$ by $r(x)=d(p,x)$, and the forward geodesic ball centered at $p$ with radius $R$ by $B(p,2R0)$. Then we obtain that

\begin{thm}\label{thm-noncpt-reduced}
	Let $(M,F,\mu)$ be a complete noncompact Finsler metric measure space without boundary, whose mixed weighted Ricci curvature satisfies that $^mRic^{N}_{\nabla r}\geq-K(2R)$ in $B(p,2R)$ with $N>n$ and $K(2R)\geq0$. Moreover, suppose the misaligment $\alpha$ of $M$ satisfies $\alpha\leq A(2R)$ in $B(p,2R)$, and suppose the non-Riemannian curvatures $U$, $\mathcal{T}$ and $\mathrm{div} C(V)=C^{ij}_{\,\,\,k|i}(V)V^k\frac{\delta}{\delta x^j}$ satisfy the norm bounds by $F(U)+F(\mathcal{T})+F^*(\mathrm{div} C(V))\leq K_0.$ Let $u$ be a bounded positive solution to the Finslerian Allen-Cahn equation \eqref{equ-FAC} on $M$ with $u\leq C$, where $C$ is a positive constant.
	
	(1) If $C\leq1$, then we have
	\begin{eqnarray*}
		\begin{split}
			\frac{F^2(\nabla u)}{u^2}+\frac{2}{3}(1-u^2)\leq (N+n)\left(\frac{C_1(1+R+R\sqrt{K(2R)})}{R^2}+2K(2R)\right),
		\end{split}
	\end{eqnarray*} 
	in $B(p,R)$, where $C_0,C_1,C_2$ are positive constants related to $A, K_0, N$.
	
	(2) If $C>1$, then we have
	\begin{eqnarray*}
		\begin{split}
			\frac{F^2(\nabla u)}{u^2}+2(1-u^2)\leq&\frac{27}{\sqrt2}(N+n)^{\frac32}C^2+2(N+n)K(2R)\\
			&+\frac{C_1}{R^2}\left(1+R+R\sqrt{K(2R)}\right),
		\end{split}
	\end{eqnarray*}
	in $B(p,R)$, where $C_0,C_1,C_2$ are positive constants depending on $A, K_0, N$.
\end{thm}

As an application of those estimates, we present the following Liouville type theorem.
\begin{thm}\label{thm-Liouville}
	Let $(M,F,\mu)$ be a forward complete Finsler metric measure space without boundary, satisfying one of the following cases.
	\item[(a)] The mixed weighted Ricci curvature of $M$ is nonnegative for some $N>n$. Moreover, $M$ admits finite misalignment and the non-Riemannian curvatures $U$, $\mathcal{T}$ and $\mathrm{div} C(V)=C^{ij}_{\,\,\,k|i}(V)V^k\frac{\delta}{\delta x^j}$ satisfy the norm bounds by $F(U)+F(\mathcal{T})+F^*(\mathrm{div} C(V))\leq K_0.$
	\item[(b)] $M$ is a compact Finsler $CD(0,N)$ space.\\	
	If $u$ is a solution to \eqref{equ-FAC} with $0<u\leq1$, then $u$ is equal to 1 identically on $M$.
\end{thm}

This manuscript is arranged as follows. In Section 2, we briefly review some Riemannian and non-Riemannian concepts on Finsler metric measure spaces. In Section 3, we introduce the Finslerian Allen-Cahn equation with its variation functional, and some regularity discussion of its solutions. Then we give the global gradient estimates of the bounded positive solutions in Section 4 with the $CD(-K,N)$ condition and the local gradient estimates in Section 5 with the mixed weighted Ricci curvature condition.

\section{Finsler metrics and Finsler metric measure spaces}

A Finsler metric space is a triple $(M,F,\mu)$, which indicates that a differential manifold is equipped with a Finsler metric $F$ and a measure $\mu$. Similar to the Riemann manifold, the Riemann curvatures on a Finsler metric measure space are still the key to determining the geometric, topological and analytical properties of the manifold. Suppose the local coordinates of the tangent bundle is $(x,y)$, where $x$ is the point on $M$ and $y$ is the direction on $T_xM$. A Finsler metric $F$ is a nonnegative function $F:TM\rightarrow [0,+\infty)$ obeying that
\begin{itemize}
	\item[(i)] $F$ is smooth and positive on $TM\setminus\{0\}$;
	\item[(ii)] $F$ is a positive homogenous norm, i.e., $F(x,ky)=kF(x,y)$ for any $(x,y)\in TM$ and for any $k>0$;
	\item[(iii)] $F$ is strongly pseudo-convex, namely, for any $(x,y)\in TM\setminus\{0\}$, the \emph{fundamental tensor} is a positive definite matrix defined by
	\begin{eqnarray}\label{Fiii}
		g_{ij}(x,y):=\frac12\frac{\partial F^2}{\partial y^i\partial y^j}(x,y).
	\end{eqnarray}
\end{itemize}

Unlike the Riemann metric, the Finsler metric is defined locally as the norm on the tangent space at each point, and globally as a metric on the pullback bundle, so there are a large number of non-Riemannian geometric quantities on a Finsler metric measure space.
The \emph{Cartan tensor} is defined by
$$C(X,Y,Z):=C_{ijk}X^iY^jZ^k=\frac{1}{4}\frac{\partial^3F^2(x,y)}{\partial y^i\partial y^j\partial y^k}X^iY^jZ^k,$$
for any local vector fields $X,Y,Z$. 

There is a unique almost $g$-compatible and torsion-free connection on the pull back tangent bundle $\pi^*TM$ of the Finsler manifold $(M,F)$ called the \emph{Chern connection}. It is determined by
\begin{align*}
	\nabla_XY-\nabla_YX&=[X,Y];\\
	Z(g_y(X,Y))-g_y(\nabla_ZX,Y)-&g_y(X,\nabla_ZY)=2C_y(\nabla_Zy,X,Y),
\end{align*}
for any $X,Y,Z\in TM\setminus\{0\}$, where $C_y$ is the Cartan tensor. The Chern connection coefficients is locally denoted by $\Gamma^i_{jk}(x,y)$ in natural coordinate system, which could induce the \emph{spray coefficients} as $G^i=\frac12\Gamma^i_{jk}y^jy^k$. The spray is given by 
\begin{eqnarray}
	G=y^i\frac{\delta}{\delta x^i}=y^i\frac{\partial}{\partial x^i}-2G^i\frac{\partial}{\partial y^i},
\end{eqnarray}
in which $\frac{\delta}{\delta x^i}=\frac{\partial}{\partial x^i}-N^j_i\frac{\partial}{\partial y^j}$, and the nonlinear connection coefficients are locally induced from the spray coefficients by $N^i_j=\frac{\partial G^i}{\partial y^j}$.

The Chern connection can define the \emph{Chern Riemannian curvature} $R$ and \emph{Chern non-Riemannian connection} $P$. Denote by $\Omega$ the curvature form of Chern connection, so that
$$\Omega(X,Y)Z=R(X,Y)Z+P(X,\nabla_Yy,Z),$$
for any $X,Y,Z\in TM\setminus\{0\}$, where locally
$$R_{j\,\,kl}^{\,\,i}=\frac{\delta \Gamma^i_{jl}}{\delta x^k}+\frac{\delta \Gamma^i_{jk}}{\delta x^l}+\Gamma^i_{km}\Gamma^{m}_{jl}-\Gamma^i_{lm}\Gamma^{m}_{jk},\quad P_{j\,\,kl}^{\,\,i}=-\frac{\partial \Gamma^i_{jk}}{\partial y^l}.$$
Customarily, we denote the horizontal Chern derivative by $``|"$ and the vertical Chern derivative by $``;"$. For example, 
$$v_{i|j}=\frac{\delta}{\delta x^j}v_i-\Gamma^k_{ij}v_k,\quad v_{i;j}=\frac{\partial }{\partial y^j}v_i,$$
for any 1-form $v=v_idx^i$ on the pull-back bundle.

The \emph{angular metric form} $h_y$ is defined by
\begin{eqnarray}\label{Def-amf}
	h_y(u, v)=g_y(u,v)-\frac{1}{F^2(y)}g_y(y, u)g_y(y, v),
\end{eqnarray}
for any $y, u, v\in T_xM$ with $y\neq 0$. Thus,
for any two linearly independent vector $y,u\in T_xM\setminus\{0\}$, which span a tangent plane $\Pi_y=\mathrm{span}\{y,u\}$, the \emph{flag curvature} with pole $y$ is defined by
\begin{eqnarray}
	K(P,y) = K(y,u):=\frac{R_y(y, u, u, y)}{F^2(y)h_y(u, u)},
\end{eqnarray}
which is locally expressed by
$$K(y,u)=\frac{-R_{ijkl}(y)y^iu^jy^ku^l}{(g_{ik}(y)g_{jl}(y)-g_{il}(y)g_{jk}(y))y^iu^jy^ku^l}.$$
The \emph{Ricci curvature} is defined by
$$Ric(y):=F^2(y)\sum_{i=1}^{n-1}K(y,e_i),$$
where $e_1,\cdots,e_{n-1},\frac{y}{F(y)}$ form an orthonormal basis of $T_xM$ with respect to $g_y$. 
The \emph{Landsberg curvature} of $(M,F)$ is given by
$$L:=L^i_{jk}\partial_i\otimes dx^j\otimes dx^k, \quad L^i_{jk}=-y^jP^{\,\,i}_{j\,\,kl}.$$
By the zero homogeneity, according to the Euler lemma, it is easy to see that 
$$C_{ijk}y^i=L_{ijk}y^i=0.$$
More Riemannian or non-Riemannian curvatures with their geometric properties can be found in \cite{BCS2000} or \cite{ShenShen2016}. 

For any given non-vanishing vector field $V$, $g_{ij}(x,V)$ could induce a Riemannian structure $g_V$ on $M$ via
\begin{eqnarray*}
	g_V(X,Y)=g_{ij}(x,V)X^iY^j,
\end{eqnarray*}
for any vector pairs $X,Y\in T_xM$.
The norm of $V$ is measured as $F^2(V)=g_V(V,V)$. Further, such an induced Riemannian metric can be checked to assign the given direction to the connection. From the asymmetry of the Finsler metric, connections in different directions at a point may not be coincident.
Noticing that the difference of two connections is a tensor. The \emph{$T$ curvature} (or \emph{tangent curvature}) introduced in Section 10.1 of \cite{ShenLec} is defined by
\begin{eqnarray}
	T_y(v):=g_y(D_vV,y)-\hat g_x(\hat D_vV,y),
\end{eqnarray}
where $v\in T_xM$, $V$ is a vector field with $V(x)=v$, and $\hat D$ denotes the Levi-Civita connection of the induced Riemannian metric $\hat g=g_Y$. The $T$ curvature vanishes if and only if the Chern non-Riemannian curvature $P$ vanishes. Similarly, let $\{e_i\}$ be an orthonormal basis with respect to metric $g(x,V)$ at point $x$, where $V$ is a fixed reference vector field. Moreover, Let $\{E_i\}$ be the local vector fields obtained by moving $\{e_i\}$ in parallel in a neighborhood of $x$ on $M$. \cite{Shen2024} defined tensor $U$ as 
\begin{eqnarray}
	U_y(W)=g(x,W)(U(y,W),W),
\end{eqnarray} 
for any local vector field $W$, with 
\begin{eqnarray}\label{def-UyW}
	U(y,W)=\sum_{i=1}^n (D^W_{e_i}E_i-\hat D_{e_i}E_i)
\end{eqnarray}
being a vector field on the sphere bundle, where $\hat D$ denotes the Levi-Civita connection of the induced Riemannian metric $\hat g=g_Y$, with $Y$ denoting the geodesic extension field of $y$ in a neighborhood of $x$. \eqref{def-UyW} is also a curvature, which may be considered as a kind of trace of the $T$ curvature.\\

Motivated by the generalization of the topological sphere theorem in Finsler geometry, H. B. Rademacher introduced the concepts of reversibility and reversible manifolds  \cite{Red2004}, which are also closely related to the analytical assumptions on Finsler manifolds.
A Finsler metric is defined to be \emph{reversible} if $F(x,V)=\bar F(x,V)$ for any point $x$ and any vector field $V$, where $\bar F(x,V):=F(x,-V)$ is called the \emph{reversed Finsler metric} of $F$. We define the \emph{reversibility} of $(M,F)$ by
\begin{eqnarray*}
	\rho:=\sup_{x\in M}\sup_{V\neq 0}\frac{F(x,V)}{\bar F(x,V)}.
\end{eqnarray*}
Obviously, $F$ is \emph{reversible} if and only if $\rho\equiv 1$. A Finsler manifold $(M,F)$ is said to have \emph{finite reversibility} if $\rho<+\infty$.

Later, K. Ball, E. Carlen and E. Lieb introduced the uniform smoothness and the uniform convexity in Banach space theory \cite{BCL1994Invent}, whose geometric explanation in Finsler geometry was given by S. Ohta \cite{Ohta2017}. We say $F$ satisfies \emph{uniform convexity} and \emph{uniform smoothness} if there exist uniform positive constants $\kappa^*$ and $\kappa$, called the \emph{uniform convexity constant} and \emph{uniform smoothness constant}, respectively, such that for any $x\in M$, $V\in T_xM\setminus\{0\}$ and $y\in T_xM$, we have 
\begin{eqnarray}\label{kappa}
	\kappa^*F^2(x,y)\leq g_V(y,y)\leq\kappa F^2(x,y),
\end{eqnarray}
where $g_V=(g_{ij}(x,V))$ is the Riemannian metric on $M$ induced from $F$ with respect to the reference vector $V$. In this situation, the reversibility $\rho$ could be controlled by $\kappa$ and $\kappa^*$ as
\begin{eqnarray}
	1\leq \rho\leq \min\{\sqrt{\kappa},\sqrt{1/\kappa^*}\}.
\end{eqnarray}
$F$ is Riemannian if and only if $\kappa=1$, if and only if $\kappa^*=1$ \cite{Ohta2017}.

The Riemannian structure is inconsistent when the reference vectors are different. For example, given three different local non-vanishing vector fields around $x$, namely, $V,W,Y$, the norm of $Y$ about $g_V$ and $g_W$ maybe not the same in general case. The ratio $g_V(Y,Y)/g_W(Y,Y)$ is a function about $V,W,Y$. 
Based on this fact, \cite{Shen2024} defined an important constant on a Finsler manifold, called the misalignment.

\begin{defn}[\cite{Shen2024}]
	For a Finsler manifold $(M,F)$, the misalignment of a Finsler metric at point $x$ is defined by
	\begin{eqnarray}
		\alpha(x)=\sup_{V,W,Y\in S_xM}\frac{g_V(Y,Y)}{g_W(Y,Y)}.
	\end{eqnarray}
	Moreover, the global misalignment of the Finsler metric is defined by
	\begin{eqnarray}\label{def-alpha}
		\alpha=\sup_{x\in M}\alpha_M(x)=\sup_{x\in M}\sup_{V,W,Y\in S_xM}\frac{g_V(Y,Y)}{g_W(Y,Y)}.
	\end{eqnarray}
\end{defn}
\cite{Shen2024} also provided some characterizations of the misalignment. Specially, a Finsler manifold $(M,F)$ is a Riemannian manifold if and only if $\alpha_M=1$. Moreover, a Finsler manifold $(M,F)$ is uniform convexity and uniform smoothness if and only if it satisfies finite misalignment.

Since that, we give an important class of Finsler manifold as the following.
\begin{defn}
	We call a Finsler manifold $(M,F)$ has finite misalignment if there is a positive constant $A$ such that $\alpha\leq A$, and has locally finite misalignment if for any compact subset $\Omega\subset M$, there is a constant $A(\Omega)$ depending on $\Omega$ such that $\alpha(x)\leq A(\Omega)$ for any $x\in \Omega$.
\end{defn}

So far, we have briefly introduced some Riemannian or non-Riemannian local quantities and tensors in Finsler geometry corresponding to Riemannian geometric quantities. Next, we will introduce more tensors related to the measure $\mu$.

For any smooth function $f:M\rightarrow \mathbb{R}$, $df$ denotes its differential 1-form and its \emph{gradient} $\nabla f$ is defined as the dual of the 1-form via the Legendre transformation, namely, $\nabla f(x):=l^{-1}(df(x))\in T_xM$. Locally it can be written as
\begin{eqnarray*}
	\nabla f=g^{ij}(x,\nabla f)\frac{\partial f}{\partial x^i}\frac{\partial }{\partial x^j}
\end{eqnarray*}
on $M_f:=\{df\neq 0\}$. The Hessian of $f$ is defined via the Chern connection by
\begin{eqnarray*}
	\nabla^2f(X,Y)=g_{\nabla f}(\nabla_X^{\nabla f}\nabla f,Y).
\end{eqnarray*}
It can be shown that $\nabla^2f(X,Y)$ is symmetric \cite{Ohta2014}.\\

For any two points $p,q$ on $M$, the \emph{distance function} is defined by
$$d_p(q):=d(p,q):=\inf_{\gamma}\int_0^1F(\gamma(t),\dot\gamma(t))dt,$$
where the infimum is taken over all the $C^1$ curves $\gamma:[0,1]\rightarrow M$ such that $\gamma(0)=p$ and $\gamma(1)=q$. 
Fixing a base point $p$ on $M$, we denote the forward distance function by $r$. That is, $r(x)=d(p,x)$, with $d$ denotes the forward distance. 

The \emph{forward distance function} $r$ is a function defined on the Finsler manifold $M$. $dr$ is a 1-form on $M$, whose dual is a gradient vector field, noted by $\nabla r$. Precisely, $\nabla r=g^{ij}(x,\nabla r)\frac{\partial r}{\partial x^i}\frac{\partial}{\partial x^j}$. Taking the Chern horizontal derivative of $\nabla r$ yields the\emph{ Hessian of distance function} $\nabla^2r$. Locally, in natural coordinate system 
\begin{eqnarray}\label{Hessianr}
	\nabla^2 r(\frac{\partial}{\partial x^i},\frac{\partial }{\partial x^j})=\nabla_j\nabla_ir=\frac{\partial^2 r}{\partial x^i\partial x^j}-\Gamma^k_{ij}(x,\nabla r)\frac{\partial r}{\partial x^k}.
\end{eqnarray} 
In \eqref{Hessianr}, the derivative is taken in the direction $\nabla r$ naturally. Generally, we can take the derivative in any direction. Suppose $V$ is a local vector field around $x$ on $M$. 



Note that the distance function may not be symmetric about $p$ and $q$ unless $F$ is reversible. A $C^2$ curve $\gamma$ is called a \emph{geodesic} if locally
$$\ddot\gamma(t)+2G^i(\gamma(t),\dot \gamma(t))=0,$$
where $G^i(x,y)$ are the \emph{spray coefficients}.
A \emph{forward geodesic ball} centered at $p$ with radius $R$ can be represented by
\begin{eqnarray*}
	B_R^+(p):=\{q\in M\,:\, d(p,q)<R\}.
\end{eqnarray*}
Adopting the exponential map, a Finsler manifold $(M,F)$ is said to be \emph{forward complete} or forward geodesically complete if the exponential map is defined on the entire $TM$. Thus, any two points in a forward complete manifold $M$ can be connected by a minimal forward geodesic. Moreover, the forward closed balls $\overline{B_R^+(p)}$ are compact.\\

A \emph{Finsler metric measure space} $(M,F,\mu)$ is a Finsler manifold equipped with an given measure $\mu$. 
In local coordinates $\{x^i\}_{i=1}^n$, we can express the volume form as $d\mu=\sigma(x)dx^1\wedge\cdots\wedge dx^n$ with some positive function $\sigma(x)$. For any $y\in T_xM\setminus\{0\}$, define 
$$\tau(x,y):=\log\frac{\sqrt{\det g_{ij}(x,y)}}{\sigma(x)},$$
which is called the \emph{distortion} of $(M,F,\mu)$. The definition of the \emph{S-curvature} is given in the following.
\begin{defn}[\cite{Shen1997}\cite{ShenShen2016}]\label{def-S}
	Suppose $(M,F,\mu)$ is a Finsler metric measure space. For any point $x\in M$, let $\gamma=\gamma(t)$ be a forward geodesic from $x$ with the initial tangent vector $\dot\gamma(0)=y$. The S-curvature of $(M,F,\mu)$ is
	\begin{eqnarray*}
		S(x,y):=\frac{d}{dt}\tau=\frac{d}{dt}(\frac12\log\det(g_{ij})-\log\sigma(x))(\gamma(t),\dot\gamma(t))|_{t=0}.
	\end{eqnarray*}
\end{defn}
Definition \ref{def-S} means that the S-curvature is the changing of distortion along the geodesic in direction $y$. Modeled on the definition of $T$ curvature in \cite{ShenLec}, \cite{Shen2024} defined the difference of $\nabla\tau$ on the tangent sphere, denoted by $\mathcal{T}$, 
\begin{defn}[\cite{Shen2024}]
	The difference of $\nabla\tau$ on the tangent bundle is a tensor denoted by $\mathcal{T}$, which is given by 
	\begin{eqnarray}\label{def-TVW}
		\mathcal{T}(V,W):=\nabla^V\tau(V)-\nabla^W\tau(W),
	\end{eqnarray}
	for vector fields $V,W$ on $M$. Locally, it is $\mathcal{T}(V,W)=\mathcal{T}_i(V,W)d x^i$, with 
	\begin{eqnarray}
		\mathcal{T}_i=\frac{\delta}{\delta x^i}\tau_k(V)-\frac{\delta}{\delta x^i}\tau_k(W).
	\end{eqnarray}
\end{defn}
Obviously, $\mathcal{T}(V,W)$ is anti-symmetric about $V$ and $W$, that is, $\mathcal{T}(V,W)=-\mathcal{T}(W,V)$ for any nonvanishing $V$ and $W$.

If in local coordinates $\{x^i\}_{i=1}^n$, expressing $d\mu=e^{\Phi}dx^1\cdots dx^n$, the \emph{divergence} of a smooth vector field $V$ can be written as
\begin{eqnarray}
	div_{\mu}V=\sum_{i=1}^n(\frac{\partial V^i}{\partial x^i}+V^i\frac{\partial \Phi}{\partial x^i}).
\end{eqnarray}
The \emph{Finsler Laplacian} of a function $f$ on $M$ could now be given by
\begin{eqnarray}\label{def-WLap}
	\Delta_{\mu}f:=div_{\mu}(\nabla f).
\end{eqnarray}
Noticing that $\Delta_{\mu}f=\Delta_{\mu}^{\nabla f}f$, where $\Delta_{\mu}^{\nabla f}f:=div_{\mu}(\nabla^{\nabla f}f)$ is in the view of weighted Laplacian with
\begin{eqnarray}
	\nabla^{\nabla f}f:=\begin{cases}
		&g^{ij}(x,\nabla f)\frac{\partial f}{\partial x^i}\frac{\partial}{\partial x^j} \quad \mbox{for } x\in M_f;\\
		&0\quad\quad\quad\quad\quad\quad\quad\,\,\, \mbox{for }x\notin M_f.
	\end{cases}
\end{eqnarray}

A Finsler Laplacian is better to be viewed in a weak sense due to the lack of regularity. Concretely, assuming $f\in W^{1,p}(M)$,
\begin{eqnarray*}
	\int_M\phi\Delta_{\mu}fd\mu=-\int_Md\phi(\nabla f)d\mu,
\end{eqnarray*}
for any test function $\phi\in C^{\infty}_0(M)$.

On the other hand, the Laplacian of a function $f$ on a Riemannian manifold is the trace of the Hessian of $f$ with respect to the Riemannian metric $g$. On a Finsler metric measure space $(M,F,\mu)$, the \emph{weighted Hessian} $\tilde{H}(f)$ of a function $f$ on $M_f=\{x\in M:\nabla f|_x \neq 0\}$ is defined in \cite{Wu2015} by
\begin{eqnarray}\label{Def-WHf}
	\tilde{H}(f)(X,Y)=\nabla^2f(X,Y)-\frac{S(\nabla f)}{n-1}h_{\nabla f}(X,Y),
\end{eqnarray}
where $h_{\nabla f}$ is the angular metric form in the direction $\nabla f$, given in \eqref{Def-amf} 
It is clear that $\tilde{H}(f)$ is still a symmetric bilinear form with
\begin{eqnarray}
	\Delta f=\mathrm{tr}_{\nabla f}\tilde{H}(f).
\end{eqnarray}

Inspired by this, \cite{Shen2024} defined the \emph{mixed weighted Hessian} \begin{eqnarray}\label{Def-MWH}
	\tilde H^V(f)(X,Y)=\nabla^2(f)(X,X)-\frac{S(\nabla f)}{n-1}h_{V,\nabla f}(X,X),
\end{eqnarray}
where $h_{V,\nabla f}$ is the \emph{mixed angular metric form in the directions $V$ and $\nabla f$}, which is defined by
\begin{eqnarray}
	h_{V,\nabla f}(X,Y)=g_V(X,Y)-\frac{1}{F_V^2(\nabla f)}g_V(X,\nabla f)g_V(Y,\nabla f),
\end{eqnarray}
for any vector $X,Y$.

It is necessary to remark that $h_{\nabla f,\nabla f}=h_{\nabla f}$, so that $\tilde H^{\nabla f}(f)=\tilde H(f)$ for any function $f$ on $M$.\\

With the assistance of the S-curvature, one can present the definition of the\emph{ weighted Ricci curvature} as the following.	
\begin{defn}[\cite{Ohta2014}\cite{ShenShen2016}]
	Given a unit vector $V\in T_xM$ and an positive number $k$, the weighted Ricci curvature is defined by
	\begin{eqnarray}
		Ric_k(V):=\begin{cases}
			&Ric(x,V)+\dot{S}(x,V)\quad\mbox{if } S(x,V)=0 \mbox{ and } k=n \mbox{ or if }k=\infty;\\
			&-\infty\quad\quad\quad\quad\quad\quad\quad\quad\quad\quad\quad\quad\mbox{if }S(x,V)\neq0 \mbox{ and if } k=n;\\
			&Ric(x,V)+\dot{S}(x,V)-\frac{S^2(x,V)}{k-n}\quad\quad\quad\quad\quad\quad\quad\mbox{ if }n<k<\infty,
		\end{cases}
	\end{eqnarray}
	where the derivative is taken along the geodesic started from $x$ in the direction of $V$.
\end{defn}
According to the definition of weighted Ricci curvature, B. Wu defined the \emph{weighted flag curvature} when $k=N\in (1,n)\cup(n,\infty)$ in \cite{Wu2015}. We have completely introduced this concept for any $k$ in \cite{Shen2024}.
\begin{defn}[\cite{Shen2024}]\label{Def-wfc}
	Let $(M,F,\mu)$ be a Finsler metric measure space,
	and $V,W\in T_xM$ be linearly independent vectors. The
	weighted flag curvature $K^k(V; W)$ is defined by
	\begin{eqnarray}
		K^k(V; W):=\begin{cases}
			&K(V; W)+\frac{\dot S(V)}{(n-1)F^2(V)}\,\,\,\mbox{if } S(x,V)=0 \mbox{ and } k=n \mbox{ or if }k=\infty;\\
			&-\infty\quad\quad\quad\quad\quad\quad\quad\quad\quad\quad\quad\,\,\mbox{if }S(x,V)\neq0 \mbox{ and if } k=n;\\
			& K(V; W)+\frac{\dot S(V)}{(n-1)F^2(V)}-\frac{S^2(V)}{(n-1)(k-n)F^2(V)}
			
			\quad\quad\quad\mbox{ if }n<k<\infty,
		\end{cases}
	\end{eqnarray}
	where the derivative is taken along the geodesic started from $x$ in the direction of $V$.
\end{defn} 
Moreover, the \emph{mixed weighted Ricci curvature} with respect to vector field $W$ has also been defined in \cite{Shen2024}, denoted by $^mRic^k_W$.
\begin{defn}[\cite{Shen2024}]\label{Def-mwrc}
	Given two unit vectors $V,W\in T_xM$ and a positive number $k$, the mixed weighted Ricci curvature $^mRic^k(V,W)=\,^mRic^k_W(V)$ is defined by
	\begin{eqnarray}
		^mRic^k_{W}(V):=\begin{cases}
			&\mathrm{tr}_WR_{V}(V)+\dot{S}(x,V)\quad\mbox{if } S(x,V)=0 \mbox{ and } k=n \mbox{ or if }k=\infty;\\
			&-\infty\quad\quad\quad\quad\quad\quad\quad\quad\quad\quad\quad\quad\mbox{if }S(x,V)\neq0 \mbox{ and if } k=n;\\
			&\mathrm{tr}_WR_{V}(V)+\dot{S}(x,V)-\frac{S^2(x,V)}{k-n}\quad\quad\quad\quad\quad\quad\quad\mbox{ if }n<k<\infty,
		\end{cases}
	\end{eqnarray}
	where the derivative is taken along the geodesic started from $x$ in the direction of $V$, and $\mathrm{tr}_WR_{V}(V)=g^{ij}(W)g_{V}(R_V(e_i,V)V,e_j)$ means taking trace of the flag curvature 
	with respect to $g(x,W)$.
\end{defn}

\begin{rem}
	\begin{itemize}
	\item[i)] The weighted Ricci curvature is a special case of the mixed weighted Ricci curvature, i.e., $Ric^k(V)=\,^mRic^k_V(V)$.
	\item[ii)] We call the mixed weighted Ricci curvature bounded from below by a constant $K$, denoted by $^mRic^k_W\geq K$, if the inequality holds that $^mRic^k_W(V)\geq KF^2(V)$.
	\item[iii)] We define the mixed weighted Ricci curvature $^mRic^k_{W}(V)$ to be the weighted Ricci curvature $Ric^k(V)$ when the reference vector $W=0$.
	\end{itemize} 
\end{rem}

Defining the function $\mathfrak{ct}_c(r)$ as
\begin{eqnarray}
	\mathfrak{ct}_c(r)=\begin{cases}
		\sqrt{c}\cot\sqrt{c}r,\quad\,\,\,\quad c>0,\\
		1/r, \quad\quad\quad\quad\quad\,\,\,\,  c=0,\\
		\sqrt{-c}\coth\sqrt{-c}r, \quad c<0.
	\end{cases}
\end{eqnarray}
the following weighted Hessian comparison theorem is cited from Theorem 3.3 in \cite{Wu2015}.
\begin{thm}[\cite{Wu2015}]\label{thm-HessComp}
	Let $(M,F,\mu)$ be a Finsler metric measure space and $r=d_F(p,\cdot)$ be the forward distance from $p$. Suppose that for
	some $N>n$, the weighted flag curvature of $M$ satisfies $K_N\geq\frac{N-1}{n-1}c$, where $c$ is a positive constant. Then, for
	any vector field $X$ on $M$, the following inequality is valid whenever $r$ is smooth. Particularly, $r < \pi/\sqrt{c}$ when $c>0$.
	\begin{eqnarray}
		\tilde H(r)(X,X)\leq \frac{N-1}{n-1} \mathfrak{ct}_c(r)(g_{\nabla r}(X,X)-g_{\nabla r}(\nabla r,X)^2).
	\end{eqnarray}
\end{thm}

Furthermore, \cite{Shen2024} obtained the following significant Laplacian comparison theorem.
\begin{thm}[\cite{Shen2024}]\label{thm-LapComp-1}
	Let $(M,F,\mu)$ be a forward complete $n$-dimensional Finsler metric measure space with finite msalignment $\alpha$.
	Denote the forward distance function by $r$ and by $V$ a fixed vector field on $M$. Suppose
	the mixed weighted Ricci curvature $^mRic^N_{\nabla r}$ of $M$ is bounded from below by $-K$ with $K>0$, for some $N>n$, as well as the non-Riemannian curvatures $U$, $\mathcal{T}$ and $\mathrm{div} C(V)=C^{ij}_{\,\,\,k|i}(V)V^k\frac{\delta}{\delta x^j}$ satisfy the norm bounds by $F(U)+F(\mathcal{T})+F^*(\mathrm{div} C(V))\leq K_0.$
	Then, by setting $l=K/C(N,\alpha)$ with $C(N,\alpha)=N+(\alpha-1)n-\alpha$, wherever $r$ is $C^2$, the nonlinear Laplacian of $r$ with reference vector $V$ satisfies 
	\begin{eqnarray}
		\Delta^Vr\leq C(N,\alpha)\mathfrak{ct}_{-l}(r)+\sqrt{\alpha}K_0.
	\end{eqnarray}
\end{thm}

\section{Finslerian Allen-Cahn equation}

We consider the Finslerian Allen-Cahn equation \eqref{equ-FAC} on a Finsler metric measure space $(M,F,\mu)$, where $\Delta^{\nabla u}$ is given in \eqref{def-WLap}. The solution to \eqref{equ-FAC} is a function $u(x)$ on $M$.

Firstly, \eqref{equ-FAC} is a critical point of an entropy functional. That is,
\begin{thm}
	Finslerian Allen-Cahn equation \eqref{equ-FAC} is the Euler-Lagrange equation of the functional
	\begin{eqnarray}
		\int_U\left(\frac12F^2(\nabla u)+W(u)\right)d\mu,\quad \mbox{with} \quad W(u)=\frac14(1-u^2)^2,
	\end{eqnarray}
	for any local regin $U\subset M$.
\end{thm}
\begin{proof}
	Choose a local function $v$ with compact support in $U$. A direction variational calculation shows that
	\begin{eqnarray}
		\begin{split}
			&\frac{\partial}{\partial t}\int_M\left[\frac12F^2(\nabla u+t\nabla v)+W(u+tv)\right]d\mu\mid_{t=0}\\
			=&\int_M\left[dv(\nabla u)-C_{\nabla u}(\nabla u,\nabla u,\nabla v)+W'(u)v\right]d\mu\\
			=&\int_U[dv(\nabla u)+vu(1-u^2)]d\mu.
		\end{split}
	\end{eqnarray}
	The conclusion holds for the arbitrariness of $v$. 
\end{proof}

We also use $\nabla$ to denote the horizontal Chern connection with respect to the direction $\nabla u$ and $\Delta=\Delta^{\nabla u}$ for short, where $u$ is the function that $\nabla$ or $\Delta$ acts on. 
Here $\nabla u$ is the gradient of $u$. 
However, for any function $u$ on manifold $M$, the two concepts are the same, i.e.,
$$\nabla_i u(x)=\frac{\delta u(x)}{\delta x^i}=\frac{\partial u}{\partial x^i}=u_i=g_{\nabla u}(\nabla u,\frac{\partial}{\partial x^i}),$$ 
where $\nabla_i u$ means the horizontal Chern derivative along $x^i$ and $\nabla u$ means the gradient of $u$.

Furthermore, the Hessian of $u$ satisfies that 
$$\nabla^2u(X,Y)=g_{\nabla u}(D^{\nabla u}_X\nabla u,Y),$$
for any $X,Y\in TM$, where $D^{\nabla u}$ is the covariant differentiation with respect to the horizontal Chern connection, whose reference vector is $\nabla u$. By the definition, it is easy to check that $\nabla^2u$ is a symmetric 2-form.

Let $H^1(M):=W^{1,2}(M)$ and $H^1_0(M)$ be the closure of $C^{\infty}_0(M)$ under the norm 
$$\|u\|_{H^1}:=\|u\|_{L^2(M)}+\frac12\|F(\nabla u)\|_{L^2(M)}+\frac12\|\overleftarrow{F}(\overleftarrow{\nabla} u)\|_{L^2(M)},$$
where $\overleftarrow{F}$ is the \emph{reverse Finsler metric}, defined by $\overleftarrow{F}(x,y):=F(x,-y)$ for all $(x,y)\in TM$, and $\overleftarrow{\nabla} u$ is the gradient of $u$ with respect to the reverse metric $\overleftarrow{F}$. In fact, $\overleftarrow{F}(\overleftarrow{\nabla} u)=F(\nabla(-u))$. Then $H^1(M)$ is a Banach space with respect to the norm $\|\cdot\|_{H^1}$.

Analog to \cite{Ohta2009-2}, one can define that
\begin{defn}
	A function $u$ on $M$ is a global solution to the Finslerian Allen-Cahn equation \eqref{equ-FAC} if $u\in H^1_0(M)$, and if, for any test function $\varphi\in H^1_0(M)$ (or $\varphi\in C^{\infty}_0(M)$), it holds that
	$$\int_M\varphi(1-u^2)ud\mu=\int_{M}d\varphi(\nabla u)d\mu.$$
\end{defn}

The following interior regularity of the solution on $M$ is followed from \cite{Ohta2009-2} and \cite{Shen2024}.

\begin{thm}
	Suppose $(M,F,\mu)$ is a Finsler metric measure space with finite reversibility $\rho<\infty$. Then one can take the continuous version of a global solution $u$ of the Finslerian Allen-Cahn equation \eqref{equ-FAC}, and it enjoys the $H_{loc}^2$-regularity and the $C^{1,\beta}$-regularity in $x$  for some $0<\beta<1$.  The elliptic regularity shows that $u$ is $C^{\infty}$ on $M_u$.\\	
\end{thm}

Given an open subset $\Omega\subset M$, a local solution to the Finslerian Allen-Cahn equation \eqref{equ-FAC} is defined by

\begin{defn}
	A function $u$ on any open subset $\Omega\subset M$ is a local solution to the Finslerian Allen-Cahn equation \eqref{equ-FAC}, if $u\in L^2_{loc}(\Omega)$ with $F^*(du)\in L^2_{loc}(\Omega)$, and if, for any local test function $\varphi\in H^1_0(\Omega)$ (or $\varphi\in C^{\infty}_0(\Omega)$), it holds that
	$$\int_{\Omega}\varphi(1-u^2)ud\mu=\int_{\Omega}d\varphi(\nabla u)d\mu.$$
\end{defn}

Let $B_{2R}:=B^+_{2R}(x_0)$ be a forward geodesic ball with radius $2R$ in $M$ centered at $x_0$, and $I$ be an open interval in $\mathbb{R}$. As the global solution, if $u$ is a local solution to \eqref{equ-FAC} on $B_{2R}$, then $u\in H^2(B_{2R})\cap C^{1,\beta}(B_{2R})$ 
with $\Delta u\in H^1(B_{2R})\cap C(B_{2R})$.

\begin{rem}
	$u$ is a solution to $\Delta u+u(1-u^2)=0$ in the distributional sense if and only if $-u$ is a solution to $\overleftarrow{\Delta}u+u(1-u^2)=0$ in the distributional sense, where $\overleftarrow{\Delta}$ is the Finsler Laplacian associated to $\overleftarrow{F}$. Since the metric $\overleftarrow{F}$ is also a Finsler metric on $M$, so the results we obtained in this manuscript are also valid for the solutions of $\overleftarrow{\Delta}u+u(1-u^2)=0$.
\end{rem}

\section{Global gradient estimates on compact Finsler manifolds}

In this section, we will show the global gradient estimates of positive solutions to the Finslerian Allen-Cahn equation \eqref{equ-FAC} on compact Finsler manifolds, with the weighted Ricci curvature bounded from below. This curvature condition is utilized widely in Finsler geometric analysis.

Suppose $u$ is a positive solution to \eqref{equ-FAC}. Consider the function $w(x)=u^{-q}$ with a constant $q$ to be determined. We now have the following lemma.

\begin{lem}
	The function $w$ satisfies that
	\begin{eqnarray}\label{Lapw-2}
		\Delta^{\nabla u} w=\frac{q+1}{q}\frac{F_{\nabla u}^2(\nabla^{\nabla u}w)}{w}+qw-qw^{\frac{q-2}{q}},
	\end{eqnarray}
	on $M_u$.
\end{lem}

\begin{proof}
	The gradient of $w$ with respect to $\nabla u$ is the dual of $dw$ by the pull-back metric $g_{\nabla u}$, namely, $\nabla^{\nabla u} w=-qu^{-q-1}\nabla u$ on $M$, and its norm with respect to $g_{\nabla u}$ is $F_{\nabla u}^2(\nabla^{\nabla u} w)=q^2u^{-2q-2}F^2(\nabla u)$. Therefore, we have that 
	\begin{eqnarray}\label{nabw/w}
		\frac {F_{\nabla u}^2(\nabla^{\nabla u} w)}{w^2}=q^2u^{-2}F^2(\nabla u).
	\end{eqnarray}
	On the other hand,
	\begin{eqnarray}\label{Lapw-1}
		\begin{split}
			\Delta^{\nabla u} w&=\mathrm{div}_{\mu}(\nabla^{\nabla u} w)=\mathrm{div}_{\mu}(-qu^{-q-1}\nabla u)\\
			&=q(q+1)u^{-q-2}F^2(\nabla u)-qu^{-q-1}\Delta u,
		\end{split}
	\end{eqnarray}
	wherever $\nabla u\neq 0$. Adopting the Allen-Cahn equation and \eqref{nabw/w}, \eqref{Lapw-1} yields \eqref{Lapw-2}
\end{proof}

Now let's consider positive solutions to the Finslerian Allen-Cahn equation on $M_u$. 
By the definition of the divergence of a vector or a 1-form, we deduce that
\begin{eqnarray}\label{Lapw-0}
	\Delta^{\nabla u} w=\mathrm{div}_{\mu}(\nabla^{\nabla u} w)=e^{-\Phi}\frac{\partial}{\partial x^i}(e^{\Phi}g^{ij}(\nabla u)w_i)=\mathrm{tr}_{\nabla u}(\nabla^{\nabla u})^2 w+S(\nabla^{\nabla u} w).
\end{eqnarray}
On the othere hand,

Now we set the function $H=\frac{F_{\nabla u}^2(\nabla^{\nabla u} w)}{w^2}+\beta(1-w^{-\frac{2}{q}})$, where $\beta$ is a positive constant to be determined. Direct calculations provide that
\begin{eqnarray}
	\nabla^{\nabla u} H=\frac{\nabla^{\nabla u} F_{\nabla u}^2(\nabla^{\nabla u} w)}{w^2}-2\frac{F_{\nabla u}^2(\nabla^{\nabla u} w)\nabla^{\nabla u} w}{w^3}+\frac{2\beta}{q}w^{-\frac{q+2}{q}}\nabla^{\nabla u} w,
\end{eqnarray}
and 
\begin{eqnarray}\label{trH-1}
	\begin{split}
		\mathrm{tr}_{\nabla u}(\nabla^{\nabla u})^2H=&\frac{2}{w^2}\left[d w(\nabla^{\nabla u}(\mathrm{tr}_{\nabla u}(\nabla^{\nabla u})^2 w))+Ric(\nabla^{\nabla u}w)+\|(\nabla^{\nabla u})^2w\|_{HS(\nabla u)}^2\right]\\
		&-\frac{8}{w^3}g_{\nabla u}((\nabla^{\nabla u})^2w,\nabla^{\nabla u}w\otimes\nabla^{\nabla u}w)\\
		&-\frac{2}{w^3}F_{\nabla u}^2(\nabla^{\nabla u}w)\mathrm{tr}_{\nabla u}(\nabla^{\nabla u})^2w+\frac{6}{w^4}F_{\nabla u}^4(\nabla^{\nabla u}w)\\
		&-\frac{2\beta(q+2)}{q^2}w^{-\frac{2(q+1)}{q}}F_{\nabla u}^2(\nabla^{\nabla u}w)+\frac{2\beta}{q}w^{-\frac{q+2}{q}}\mathrm{tr}_{\nabla u}(\nabla^{\nabla u})^2w
	\end{split}
\end{eqnarray}
on $M_u$, where we have employed the Ricci-type identity in \cite{Shen2018}. On the other hand, it satisfies that 
\begin{eqnarray*}
	\nabla^{\nabla u}(S(\nabla^{\nabla u}w))=\dot S(\nabla^{\nabla u}w)+g_{\nabla u}((\nabla^{\nabla u})^2w,\nabla^{\nabla u}\tau\otimes\nabla^{\nabla u}w),
\end{eqnarray*} 
as well as
\begin{eqnarray}\label{nablaS}
	\begin{split}
		\nabla^{\nabla u}(\mathrm{tr}_{\nabla u}(\nabla^{\nabla u})^2w)=&\nabla^{\nabla u}(\Delta^{\nabla u} w)+\dot S(\nabla^{\nabla u}w)\\
		&+g_{\nabla u}((\nabla^{\nabla u})^2w,\nabla^{\nabla u}\tau\otimes\nabla^{\nabla u}w).
	\end{split}
\end{eqnarray} 

Furthermore, with the assistance of \eqref{Lapw-0} and \eqref{nablaS}, \eqref{trH-1} implies that
\begin{eqnarray}\label{trH-2}
	\begin{split}
		\mathrm{tr}_{\nabla u}(\nabla^{\nabla u})^2H=&\frac{2}{w^2}[dw(\nabla^{\nabla u}(\Delta^{\nabla u} w))+\dot S(\nabla^{\nabla u}w)+Ric(\nabla^{\nabla u}w)\\
		&\quad\,\,+g_{\nabla u}((\nabla^{\nabla u})^2w,\nabla^{\nabla u}\tau\otimes\nabla^{\nabla u}w)+\|(\nabla^{\nabla u})^2w\|^2_{HS(\nabla u)}]\\
		&+\frac{6}{w^4}F_{\nabla u}^4(\nabla^{\nabla u}w)-8\frac{g_{\nabla u}((\nabla^{\nabla u})^2w,\nabla^{\nabla u}w\otimes\nabla^{\nabla u}w)}{w^3}\\
		&-2\frac{F_{\nabla u}^2(\nabla^{\nabla u}w)(\Delta^{\nabla u} w+S(\nabla^{\nabla u}w))}{w^3}\\
		&-\frac{2\beta(q+2)}{q^2}w^{-\frac{2(q+1)}{q}}F_{\nabla u}^2(\nabla^{\nabla u}w)\\
		&+\frac{2\beta}{q}w^{-\frac{q+2}{q}}(\Delta^{\nabla u} w+S(\nabla^{\nabla u}w)),
	\end{split}
\end{eqnarray}
whereiver $\nabla u\neq 0$.

By a basic computation of Finsler geometry, we find that
\begin{eqnarray}\label{LapH-1}
	\Delta^{\nabla u}H=\mathrm{div}_{\mu}(\nabla^{\nabla u} H)=&\mathrm{tr}_{\nabla u}(\nabla^{\nabla u})^2H-d\tau(\nabla^{\nabla u} H)+2C^{\nabla u}_{\nabla^2u}(\nabla^{\nabla u }H),
\end{eqnarray}
where we have already employed the fact that $C(\nabla u,\cdot,\cdot)=0$, and $C^{\nabla u}_{\nabla^2u}(\nabla^{\nabla u}H)=u^k_{\,\,\,|i}C^{ij}_k(\nabla u)H_j$. The notation $``|"$ here means the horizontal derivative with respect to the Chern connection in the direction $\nabla u$. However, noticing that we always take the direction of $\nabla u$, thus, 
\begin{eqnarray*}
	\begin{split}
		C^{\nabla u}_{\nabla^2u}(\nabla^{\nabla u}H)&=u^k_{\,\,\,|i}C^{ij}_k(\nabla u)\left[2w^lw_{l|j}w^{-2}-2F_{\nabla u}^2(\nabla^{\nabla u}w)w^{-3}w_j+\frac{2}{q}\beta w^{-\frac{q+2}{q}}w_j\right]\\
		&=2\frac{u^k_{\,\,\,|i}}{w^2}\left[(C^{ij}_kw_j)_{|l}w^l-C^{ij}_{k|l}w_jw^l\right]\\
		&=-2\frac{u^k_{\,\,\,|i}}{w^2}L^{ij}_kw_i=0,
	\end{split}
\end{eqnarray*}
where $w^j=g^{ij}(\nabla u)w_i$. Therefore, \eqref{LapH-1} is equal to
\begin{eqnarray}\label{LapH-2}
	\Delta^{\nabla u}H=\mathrm{div}_{\mu}(\nabla^{\nabla u} H)=\mathrm{tr}_{\nabla u}(\nabla^{\nabla u})^2H-d\tau(\nabla^{\nabla u} H).
\end{eqnarray}

Moreover, it follows from \eqref{trH-2} and \eqref{LapH-2} that
\begin{eqnarray}\label{LapH-3}
	\begin{split}
		\Delta^{\nabla u}H=&\frac{2}{w^2}[g_{\nabla u}(\nabla^{\nabla u}(\Delta^{\nabla u} w),\nabla^{\nabla u}w)+Ric^{\infty}(\nabla^{\nabla u}w)+\|(\nabla^{\nabla u})^2w\|_{HS(\nabla u)}^2]\\
		&-8\frac{g_{\nabla u}((\nabla^{\nabla u})^2w,\nabla^{\nabla u}w\otimes\nabla^{\nabla u}w)}{w^3}-2\frac{F_{\nabla u}^2(\nabla^{\nabla u}w)\Delta^{\nabla u} w}{w^3}\\
		&+6\frac{F_{\nabla u}^4(\nabla^{\nabla u}w)}{w^4}-\frac{2\beta(q+2)}{q^2}w^{-\frac{2(q+1)}{q}}F_{\nabla u}^2(\nabla^{\nabla u}w)+\frac{2\beta}{q}w^{-\frac{q+2}{q}}\Delta^{\nabla u} w,
	\end{split}
\end{eqnarray}
by noticing 
\begin{eqnarray*}
	\begin{split}
		d\tau(\nabla^{\nabla u}H)=2d\tau&\left[\frac{g_{\nabla u}((\nabla^{\nabla u})^2w,\nabla^{\nabla u}w)}{w^2}\right.\\
		&\,\,\,\left.-\frac{F_{\nabla u}^2(\nabla^{\nabla u}w)}{w^3}\nabla^{\nabla u}w+\frac{\beta}{q}w^{-\frac{2+q}{q}}\nabla^{\nabla u}w\right].
	\end{split}
\end{eqnarray*}

We now estimate the terms including $\Delta^{\nabla u} w$ in \eqref{LapH-3} by adopting \eqref{Lapw-2}. That is,
\begin{eqnarray}\label{Delta-1}
	\begin{split}
		g_{\nabla u}(\nabla^{\nabla u}(\Delta^{\nabla u} w),\nabla^{\nabla u}w)=&(q-(q-2)w^{-\frac{2}{q}})\frac{F_{\nabla u}^2(\nabla^{\nabla u}w)}{w^2}-\frac{(q+1)}{q}\frac{F_{\nabla u}^4(\nabla^{\nabla u}w)}{w^4}\\
		&+\frac{2(q+1)}{q}\frac{g_{\nabla u}((\nabla^{\nabla u})^2w,\nabla^{\nabla u}w\otimes \nabla^{\nabla u}w)}{w^3},
	\end{split}
\end{eqnarray}	
\begin{eqnarray}\label{Delta-2}
	\frac{F_{\nabla u}^2(\nabla^{\nabla u}w)(\Delta^{\nabla u} w)}{w^3}=\frac{(q+1)}{q}\frac{F_{\nabla u}^4(\nabla^{\nabla u}w)}{w^4}+q(1-w^{-\frac{2}{q}})\frac{F_{\nabla u}^2(\nabla^{\nabla u}w)}{w^2},
\end{eqnarray}
as well as
\begin{eqnarray}\label{Delta-3}
	\frac{\beta}{q}w^{-\frac{q+2}{q}}\Delta^{\nabla u} w=\frac{\beta(q+1)}{q^2}\frac{F_{\nabla u}^2(\nabla^{\nabla u}w)w^{-\frac{2}{q}}}{w^2}+\beta w^{-\frac{2}{q}}-\beta w^{-\frac{4}{q}},
\end{eqnarray}
Plugging \eqref{Delta-1}-\eqref{Delta-3} into \eqref{LapH-3} yields
\begin{eqnarray}\label{LapH-4}
	\begin{split}
		\Delta^{\nabla u}H=&2\frac{Ric^{\infty}(\nabla^{\nabla u}w)+\|(\nabla^{\nabla u})^2 w\|_{HS(\nabla u)}^2}{w^2}+\left(4-\frac{2\beta}{q^2}\right)w^{-\frac{2}{q}}\frac{F_{\nabla u}^2(\nabla^{\nabla u}w)}{w^2}\\
		&+\frac{4(1-q)}{q}\frac{g_{\nabla u}((\nabla^{\nabla u})^2w,\nabla^{\nabla u}w\otimes \nabla^{\nabla u}w)}{w^3}+\frac{2(q-2)}{q}\frac{F_{\nabla u}^4(\nabla^{\nabla u}w)}{w^4}\\
		&+2\beta w^{-\frac{2}{q}}(1-w^{-\frac{2}{q}}),
	\end{split}
\end{eqnarray}
on $M_u$.

Utilizing the H\"older inequality, one could see that
\begin{eqnarray*}
	2\frac{g_{\nabla u}((\nabla^{\nabla u})^2w,\nabla^{\nabla u}w\otimes \nabla^{\nabla u}w)}{w^3}\leq\frac{\epsilon\|(\nabla^{\nabla u})^2 w\|^2_{HS(\nabla u)}}{w^2}+\frac{F_{\nabla u}^4(\nabla^{\nabla u}w)}{\epsilon w^4}.
\end{eqnarray*}
Hence
\begin{eqnarray}\label{delta-4}
	\begin{split}
		&\frac{4(1-q)}{q}\frac{g_{\nabla u}((\nabla^{\nabla u})^2w,\nabla^{\nabla u}w\otimes \nabla^{\nabla u}w)}{w^3}\\
		&+2\frac{\|(\nabla^{\nabla u})^2 w\|_{HS(\nabla u)}^2}{w^2}+\frac{2(q-2)}{q}\frac{F_{\nabla u}^4(\nabla^{\nabla u}w)}{w^4}\\
		\geq& \frac{4}{q}\frac{g_{\nabla u}((\nabla^{\nabla u})^2w,\nabla^{\nabla u}w\otimes \nabla^{\nabla u}w)}{w^3}\\
		&+2(1-\epsilon)\frac{\|(\nabla^{\nabla u})^2 w\|_{HS(\nabla u)}^2}{w^2}+2(1-\frac{1}{\epsilon}-\frac{2}{q})\frac{F_{\nabla u}^4(\nabla^{\nabla u}w)}{w^4},
	\end{split}
\end{eqnarray}
in which $0<\epsilon<1$.

According to the equality 
$$\frac{(a+b)^2}{n}=\frac{a^2}{N}-\frac{b^2}{N-n}+\frac{N(N-n)}{n}(\frac{a}{N}+\frac{b}{N-n})^2,$$
for any $N>n$, the following inequality holds by substituting $a$, $b$ and $N$ by $\Delta^{\nabla u} w$, $\dot S(\nabla^{\nabla u}w)$ and $N-\epsilon(N-n)$, respectively,
\begin{eqnarray}\label{cur-geq}
	\begin{split}
		\|(\nabla^{\nabla u})^2 w\|^2_{HS(\nabla u)}&\geq \frac{(\mathrm{tr}_{\nabla u}(\nabla^{\nabla u})^2w)^2}{n}=\frac{(\Delta^{\nabla u} w+S(\nabla^{\nabla u}w))^2}{n}\\
		&\geq\frac{(\Delta^{\nabla u} w)^2}{N-\epsilon(N-n)}-\frac{(S(\nabla^{\nabla u}w))^2}{(1-\epsilon)(N-n)}.
	\end{split}	
\end{eqnarray}
At last, it is easy to verify that
\begin{eqnarray}\label{Hw-eq}
	\begin{split}
		g_{\nabla u}(\nabla^{\nabla u} H,\nabla^{\nabla u}\log w)=2&\left(\frac{g_{\nabla u}((\nabla^{\nabla u})^2w,\nabla^{\nabla u}w\otimes \nabla^{\nabla u}w)}{w^3}\right.\\
		&\quad\left.-\frac{F_{\nabla u}^4(\nabla^{\nabla u}w)}{w^4}+\frac{\beta}{q}w^{-\frac{2}{q}}\frac{F_{\nabla u}^2(\nabla^{\nabla u}w)}{w^2}\right).
	\end{split}
\end{eqnarray}

It is deduced from \eqref{LapH-4} by combining \eqref{delta-4}-\eqref{Hw-eq} and employing the definition of weighted Ricci curvature that 
\begin{eqnarray}\label{LapH-5}
	\begin{split}
		\Delta^{\nabla u}H\geq&\frac{2(1-\epsilon)}{N-\epsilon(N-n)}\frac{(\Delta^{\nabla u} w)^2}{w^2}-2(\frac{1}{\epsilon}-1)\frac{F_{\nabla u}^4(\nabla^{\nabla u}w)}{w^4}\\
		&+\frac{2}{q}g_{\nabla u}(\nabla^{\nabla u} H,\nabla^{\nabla u}\log w)+2\frac{Ric^N(\nabla^{\nabla u}w)}{w^2}\\
		&+(4-6\frac{\beta}{q^2})\frac{F_{\nabla u}^2(\nabla^{\nabla u}w)}{w^2}w^{-\frac{2}{q}}+2\beta w^{-\frac{2}{q}}(1-w^{-\frac{2}{q}}),
	\end{split}	
\end{eqnarray}
on $M_u$.
Noticing from \eqref{Lapw-2} and the expression of $H$, we get that
\begin{eqnarray}\label{Deltaw-re}
	\frac{\Delta^{\nabla u} w}{w}=\frac{q}{\beta}H+(\frac{q+1}{q}-\frac{q}{\beta})\frac{F_{\nabla u}^2(\nabla^{\nabla u}w)}{w^2}.
\end{eqnarray}
Plugging \eqref{Deltaw-re} into \eqref{LapH-5}, and setting $\beta=sq^2$, we arrive at the following lemma.

\begin{lem}\label{lem1}
	Let $(M,F,\mu)$ be a forward complete Finsler metric measure space, and denote $M_u=\{x\in M \mid \nabla u(x)\neq 0\}$. For $H=\frac{F_{\nabla u}^2(\nabla^{\nabla u}w)}{w^2}+\beta(1-w^{-\frac{2}{q}})$ with $\beta=sq^2$, it holds on $M_u$ that
	\begin{eqnarray}\label{LapH-6}
		\begin{split}
			\Delta^{\nabla u}H\geq&\frac{2(1-\epsilon)}{N-\epsilon(N-n)}\frac{H^2}{s^2q^2}+\left[\frac{2(1-\epsilon)(sq+s-1)^2}{(N-\epsilon(N-n))s^2q^2}-2(\frac{1}{\epsilon}-1)\right]\frac{F_{\nabla u}^4(\nabla^{\nabla u}w)}{w^4}\\
			&+\frac{4(1-\epsilon)(sq+s-1)}{((1-\epsilon)N+\epsilon n)s^2q^2}H\frac{F_{\nabla u}^2(\nabla^{\nabla u}w)}{w^2}+\frac{2}{q}g_{\nabla u}(\nabla^{\nabla u} H,\nabla^{\nabla u}\log w)\\
			&+2w^{-\frac{2}{q}}H+(2-6s)\frac{F_{\nabla u}^2(\nabla^{\nabla u}w)}{w^2}w^{-\frac{2}{q}}+2\frac{Ric^N(\nabla^{\nabla u}w)}{w^2}.
		\end{split}	
	\end{eqnarray}
\end{lem}

With the assistance of Lemma \ref{lem1}, we can prove the global gradient estimates on compact Finsler metric measure spaces directly, with the weighted Ricci curvature bounded from below. Precisely, we have 
\begin{thm}\label{thm-cpt-1}
	Let $(M,F,\mu)$ be a compact Finsler metric measure space whose weighted Ricci curvature satisfies $Ric^{N}\geq -K$, with some positive constant $K$. Assume the bound of the reversibility on $M$ is $\rho_0$. Suppose $u$ is a bounded positive solution to the Finslerian Allen-Cahn equation \eqref{equ-FAC} on $M$ with $u\leq C$ where $C$ is a positive constant. 
	
	(1) If $C\leq 1$, then we have 
	$$	\frac{F^2(\nabla u)}{u^2}+\frac{2}{3}(1-u^2)\leq \frac{4\rho_0^2K[N-\epsilon(N-n)]}{3(1-\epsilon)},$$
	for any $0<\epsilon<1$.
	
	(2) If $C\geq 1$, then we have 
	$$\frac{F^2(\nabla u)}{u^2}+s(1-u^2)\leq\frac{[N-\epsilon(N-n)]}{(1-\epsilon)(1-s)^2}s^2\rho_0^2K+\frac{s}{q}\sqrt{\frac{[N-\epsilon(N-n)]}{(1-\epsilon)}}C^2,$$
	for any $0<\epsilon<1$, $s>1$ and $0\leq q\leq \min\left\{1,\frac{s-1}{s}\left[\frac{1-\epsilon}{\epsilon}+\frac{(3s-1)^2}{2}\right]^{-1}\frac{2(1-\epsilon)}{N-\epsilon(N-n)}\right\}$.
\end{thm}

\begin{proof}
	Lemma \ref{lem1} implies that \eqref{LapH-6} holds in the distributional sense on $M$, i.e., for any nonnegative test function $\phi$, $H$ satisfies that 
	\begin{eqnarray}\label{beta-int-ineq}
		-\int_Md\phi(\nabla^{\nabla u}H)d\mu\geq\int_M\phi \beta d\mu,
	\end{eqnarray}
	where $\beta$ denotes the RHS of \eqref{LapH-6}.
	
	Suppose $x_0$ is the maximum point of $H$ on $M$. Without loss of generality, one could assume $H(x_0)\geq 0$, otherwise the result will be satisfied trivially. We claim that $\beta(x_0)\leq 0$.
	
	If not, $\beta$ is strictly positive at $x_0$, so that $\beta(x)$ is positive in a small neighborhood of $x_0$ on $M$, which may be denoted by $U$. Chosen a test function $\phi$ whose compact support set is contained in $U$, we know from \eqref{beta-int-ineq} that $H$ is a weak, local subharmonic function in a neighborhood $V\subset U$. It is a contradiction because $x_0$ is a inner point of $V$.  
	
	Since $\beta(x_0)\leq 0$ and $\nabla^{\nabla u}H(x_0)=0$, we arrive at	
	\begin{eqnarray}\label{LapH-cpt-1}
		\begin{split}
			0\geq&\frac{2(1-\epsilon)}{N-\epsilon(N-n)}\frac{H^2}{s^2q^2}+\left[\frac{2(1-\epsilon)(sq+s-1)^2}{(N-\epsilon(N-n))s^2q^2}-2(\frac{1}{\epsilon}-1)\right]\frac{F_{\nabla u}^4(\nabla^{\nabla u}w)}{w^4}\\
			&+\frac{4(1-\epsilon)(sq+s-1)}{((1-\epsilon)N+\epsilon n)s^2q^2}H\frac{F_{\nabla u}^2(\nabla^{\nabla u}w)}{w^2}+2w^{-\frac{2}{q}}H\\
			&+(2-6s)\frac{F_{\nabla u}^2(\nabla^{\nabla u}w)}{w^2}w^{-\frac{2}{q}}+2\frac{Ric^N(\nabla^{\nabla u}w)}{w^2}.
		\end{split}	
	\end{eqnarray}
	
	We consider \eqref{LapH-cpt-1} in two cases: (1) $C\leq 1$ and (2) $C>1$.\\
	
	\textit{Case 1.} $u\leq 1$ in this case. It is easy to see that $w\geq1$ so that $1-w^{-\frac2q}\geq0$ and 
	\begin{eqnarray}\label{cptC1-1}
		\frac{Ric^N(\nabla^{\nabla u}w)}{w^2}\geq -K\frac{F^2(\nabla^{\nabla u}w)}{w^2}\geq -K\rho_0^2H.
	\end{eqnarray}
	Meanwhile, applying the H\"older inequality on the following terms gives that
	\begin{eqnarray}\label{cptc1-2}
		\begin{split}
			&2w^{-\frac{2}{q}}H+(2-6s)\frac{F_{\nabla u}^2(\nabla^{\nabla u}w)}{w^2}w^{-\frac{2}{q}}\\
			=& (4-6s)w^{-\frac{2}{q}}H+(6s-2)sq^2(1-w^{-\frac2q})w^{-\frac2q}\\
			\geq&-(6s-4)H,
		\end{split}
	\end{eqnarray}
	provided $s\geq \frac{1}{3}$.
	
	We may choose $s=\frac23$ so that the RHS of \eqref{cptc1-2} is equal to 0, and choose 
	$0<q\leq \frac{1}{2}\left(1+\sqrt{\frac{[N-\epsilon(N-n)]}{\epsilon}}\right)^{-1}$
	so that the second term on the RHS of \eqref{LapH-cpt-1} in this case is nonnegative. Furthermore, the first and third terms on the RHS of \eqref{LapH-cpt-1} could be bounded from below by
	\begin{eqnarray*}
		\begin{split}
			&\frac{2(1-\epsilon)}{N-\epsilon(N-n)}\frac{H^2}{s^2q^2}-\frac{4(1-\epsilon)(sq+s-1)}{(N-\epsilon(N-n))s^2q^2}H\frac{F_{\nabla u}^2(\nabla^{\nabla u}w)}{w^2}\\
			=&\frac{2(1-\epsilon)}{(N-\epsilon(N-n))s^2q^2}[H^2+2(sq+s-1)H\frac{F_{\nabla u}^2(\nabla^{\nabla u}w)}{w^2}]\\
			\geq& \frac{3(1-\epsilon)}{2(N-\epsilon(N-n))q^2}H^2. 
		\end{split}
	\end{eqnarray*}

	Plugging \eqref{cptC1-1} and \eqref{cptc1-2} into \eqref{LapH-cpt-1} with the above estimates gives that
	\begin{eqnarray*}
		0\geq \frac{3(1-\epsilon)}{2q^2(N-\epsilon(N-n))}H^2-2\rho_0^2KH,
	\end{eqnarray*}
	which is equal to
	\begin{eqnarray}\label{ineq-H-cptC1}
		H\leq \frac{4\rho_0^2Kq^2[N-\epsilon(N-n)]}{3(1-\epsilon)}.
	\end{eqnarray}
	
	\textit{Case 2.} when $C>1$, from the curvature condition and the H\"older inequality, we know that
	\begin{eqnarray}\label{cptC2-est-1}
		\begin{split}
			&\quad\frac{2Ric^N(\nabla^{\nabla u}w)}{w^2}\geq-2\rho_0^2K\frac{F_{\nabla u}^2(\nabla^{\nabla u}w)}{w^2}\\
			&\geq-\frac{2(1-\epsilon)(s-1)^2}{(N-\epsilon(N-n))s^2q^2}\frac{F_{\nabla u}^4(\nabla^{\nabla u}w)}{w^4}-\frac{(N-\epsilon(N-n))s^2q^2}{2(1-\epsilon)(s-1)^2}\rho_0^4K^2,
		\end{split}
	\end{eqnarray}
	and 
	\begin{eqnarray}\label{cptC2-est-2}
		(6s-2)\frac{F_{\nabla u}^2(\nabla^{\nabla u}w)}{w^2}w^{-\frac{2}{q}}\leq(3s-1)^2\frac{F_{\nabla u}^4(\nabla^{\nabla u}w)}{w^4}+C^4,
	\end{eqnarray}
	where we have used that $w^{-\frac{4}{q}}=u^4\leq C^4$. 
	Now we choose $s>1$ and $q>0$ (e.g. $0\leq q\leq \min\{1,\frac{s-1}{s}\left[\frac{1-\epsilon}{\epsilon}+\frac{(3s-1)^2}{2}\right]^{-1}\frac{2(1-\epsilon)}{N-\epsilon(N-n)}\}$ for any given $s>1$ and $0<\epsilon<1$) such that 
	\begin{eqnarray}\label{cptC2-est-3}
		\frac{2(1-\epsilon)}{(N-\epsilon(N-n))}\frac{s-1}{sq}\geq\frac{1}{\epsilon}-1+\frac{(3s-1)^2}{2}.
	\end{eqnarray}
	We plug \eqref{cptC2-est-1}-\eqref{cptC2-est-3} into \eqref{LapH-cpt-1} to obtain that
	\begin{eqnarray*}
		\begin{split}
			\frac{2(1-\epsilon)}{N-\epsilon(N-n)}\frac{H^2}{s^2q^2}\leq&\frac{(N-\epsilon(N-n))s^2q^2}{2(1-\epsilon)(s-1)^2}\rho_0^4K^2+C^4\\
			&-\frac{4(1-\epsilon)(sq+s-1)}{((1-\epsilon)N+\epsilon n)s^2q^2}H\frac{F_{\nabla u}^2(\nabla^{\nabla u}w)}{w^2}-2w^{-\frac{2}{q}}H\\
			\leq&\frac{(N-\epsilon(N-n))s^2q^2}{2(1-\epsilon)(s-1)^2}\rho_0^4K^2+C^4,
		\end{split}
	\end{eqnarray*}
	otherwise the conclusion holds trivially in this case. It implies 
	\begin{eqnarray}\label{ineq-H-cptC2}
		H\leq\frac{[N-\epsilon(N-n)]}{2(1-\epsilon)(s-1)}s^2q^2\rho_0^2K+\sqrt{\frac{[N-\epsilon(N-n)]}{2(1-\epsilon)}}sqC^2.
	\end{eqnarray}

	Combining Cases 1 and 2, we achieve the global gradient estimates on compact Finsler manifolds by taking $s=\frac23$ in \textit{case 1.} 
\end{proof}

\begin{rem}
	Theorem \ref{thm-cpt-reduced} follows by taking $\epsilon=\frac12$ and $q=\frac{1}{\sqrt{N+n}}$ 
	in \textit{case 1.} as well as taking $s=2$, $\epsilon=\frac12$ and $q=\frac{2}{27(N+n)}$ 
	in \textit{case 2.} of Theorem \ref{thm-cpt-1}.
\end{rem}

\section{Local gradient estimates on forward complete Finsler manifolds}

In this section, we prove the local gradient estimates on forward complete Finsler metric measure spaces with the assistance of Lemma \ref{lem1} and the Comparison theorem. 

\begin{thm}\label{thm-noncpt-1}
	Let $(M,F,\mu)$ be a complete noncompact Finsler metric measure space. Denote by $B(p,2R)$ the forward geodesic ball of radius $2R$ centered at $p$ on $M$. Suppose the mixed weighted Ricci curvature satisfies that $^mRic^{N}_{\nabla r}\geq-K(2R)$ in $B(p,2R)$ with $K(2R)\geq0$, and the misaligment $\alpha$ satisfies $\alpha\leq A(2R)$ in $B(p,2R)$. Moreover, the non-Riemannian tensors satisfy $F(U)+F^*(\mathcal{T})+F(\mathrm{div} C(V))\leq K_0$. Let $u$ be a bounded positive solution to the Finslerian Allen-Cahn equation \eqref{equ-FAC} in $B(p,2R)$ with $u\leq C$ where $C$ is a positive constant.
	
	(1) If $C\leq1$, then we have
	\begin{eqnarray*}
		\begin{split}
			\frac{F^2(\nabla u)}{u^2}&+\frac{2}{3}(1-u^2)\leq \frac{[N-\epsilon(N-n)]}{1-\epsilon}\left\{\frac{2C_1^2+C_2}{R^2}+2AK(2R)\right.\\
			&+\frac{2(N+1-\epsilon(N+1-n))}{1-\epsilon}\frac{C_1^2}{R^2}\\
			&\left.+\frac{C_1}{R}\left[\sqrt{K(2R)C(N,A)}\coth\left(R\sqrt{\frac{K(2R)}{C(N,A)}}\right)+C_0(K_0,A)\right]\right\}
		\end{split}
	\end{eqnarray*}  
	in $B(p,R)$ for any $0<\epsilon<1$, where $C_1,C_2$ are positive constants.
	
	(2) If $C>1$, then we have
	\begin{eqnarray*}
		\begin{split}
			\frac{F^2(\nabla u)}{u^2}&+s(1-u^2)\leq\frac{(N-\epsilon(N-n))s^2}{2(1-\epsilon)R^2}\left\{\frac{(N-\epsilon(N-n))s^2C_1^2}{4(1-\epsilon)(sq+s-1)}+2C_1^2+C_2\right.\\
			&\left.+\frac{AK(2R)R^2}{s-1}+C_1R\left[\sqrt{K(2R)C(N,A)}\coth\left(R\sqrt{\frac{K(2R)}{C(N,A)}}\right)+C_0(K_0,A)\right]\right\}\\
			&+\frac{s}{q}\sqrt{\frac{N-\epsilon(N-n)}{2(1-\epsilon)}}C^2,
		\end{split}
	\end{eqnarray*}
	in $B(p,R)$ for any $0<\epsilon<1$ and $0\leq q\leq \frac{2\epsilon(1-\epsilon)}{(2+23\epsilon)(N-\epsilon(N-n))}$, where $C_1,C_2$ are positive constants.
\end{thm}

\begin{proof}
	Let $\tilde\phi\in C^2[0,+\infty)$ such that $\tilde \phi(r)=1$ when $r\leq 1$ and $\tilde\phi(r)=0$ whenever $r\geq 2$. Hence $0\leq \tilde\phi(r)\leq 1$. Moreover, we may assume that
	$$-C_1\leq\frac{\tilde \phi'(r)}{\tilde \phi^{\frac{1}{2}}(r)}\leq0,\quad \tilde\phi''(r)\geq-C_2.$$
	For a fixed point $p$, denote by $r(x)$ the forward distance function from $p$ to any point $x$. We define the cut-off function by $\phi(x)=\tilde\phi(\frac{r(x)}{R})$, so that 
	$$F_{\nabla u}(\nabla^{\nabla u} \phi)\leq\frac{\sqrt\alpha C_1}{R}\phi^{\frac12}.$$
	
	Considering the function $\phi H$, 
	one could find that
	\begin{eqnarray}\label{DeltaphiH-1}
		\Delta^{\nabla u}(\phi H)=H\Delta^{\nabla u}\phi+\phi\Delta^{\nabla u}H+\frac{2}{\phi}g_{\nabla u}(\nabla^{\nabla u}(\phi H),\nabla^{\nabla u}\phi)-\frac{2H}{\phi}F_{\nabla u}^2(\nabla^{\nabla u}\phi),
	\end{eqnarray}
	where $\nabla^{\nabla u}\phi=\tilde\phi'\frac{\nabla^{\nabla u} r}{R}$, $\Delta^{\nabla u}\phi=\tilde\phi''\frac{F^2_{\nabla u}(\nabla^{\nabla u}r)}{R^2}+\tilde\phi'\frac{\Delta^{\nabla u}r}{R}$. Plugging these two terms into \eqref{DeltaphiH-1} gives that
	\begin{eqnarray}\label{DeltaphiH-2}
		\begin{split}
			\phi\Delta^{\nabla u}(\phi H)=&\phi H(\Delta^{\nabla u}\phi-\frac{2}{\phi}F_{\nabla u}^2(\nabla^{\nabla u}\phi))\\
			&+2\phi g_{\nabla u}(\nabla^{\nabla u}(\phi H),\nabla^{\nabla u}\log\phi)+\phi^2\Delta^{\nabla u}H\\
			\geq&\phi H\left[-\frac{\alpha C_2}{R^2}-\frac{2\alpha C_1^2}{R^2}-\frac{C_1}{R}(C(N,\alpha)\mathfrak{ct}_l(r)+\sqrt{\alpha}K_0)\right]\\
			&+2\phi g_{\nabla u}(\nabla^{\nabla u}(\phi H),\nabla^{\nabla u}\log\phi)+\phi^2\Delta^{\nabla u}H,
		\end{split}
	\end{eqnarray}
	on $M_u$. \eqref{DeltaphiH-2} provides by taking the estimation of $\Delta^{\nabla u}H$ in Lemma \ref{lem1} that
	\begin{eqnarray}\label{DeltaphiH-3}
		\begin{split}
			\phi\Delta^{\nabla u}(\phi H)
			\geq&2\phi g_{\nabla u}(\nabla^{\nabla u}(\phi H),\nabla^{\nabla u}\log\phi)+\frac{2\phi}{q}g_{\nabla u}(\nabla^{\nabla u}(\phi H),\nabla^{\nabla u}\log w)\\
			&-B\phi H+\frac{2(1-\epsilon)}{N-\epsilon(N-n)}\frac{(\phi H)^2}{s^2q^2}+2\phi^2\frac{Ric^N(\nabla^{\nabla u}w)}{w^2}\\
			&+\phi^2\left[\frac{2(1-\epsilon)}{N-\epsilon(N-n)}\frac{(sq+s-1)^2}{s^2q^2}-2(\frac{1}{\epsilon}-1)\right]\frac{F_{\nabla u}^4(\nabla^{\nabla u}w)}{w^4}\\
			&+\phi\frac{4(1-\epsilon)(sq+s-1)}{N-\epsilon(N-n)s^2q^2}(\phi H)\frac{F_{\nabla u}^2(\nabla^{\nabla u}w)}{w^2}+2w^{-\frac{2}{q}}\phi^2H\\
			&+(2-6s)\phi^2w^{-\frac{2}{q}}\frac{F_{\nabla u}^2(\nabla^{\nabla u}w)}{w^2}-\frac2q(\phi H)g_{\nabla u}(\nabla^{\nabla u} \phi,\nabla^{\nabla u}\log w),
		\end{split}
	\end{eqnarray}
	where $B=\frac{2\alpha C_1^2}{R^2}+\frac{C_1}{R}\left[\sqrt{K(2R)C(N,A)}\coth\left(R\sqrt{\frac{K(2R)}{C(N,A)}}\right)+C_0(K_0,A)\right]+\frac{\alpha C_2}{R^2}$.
	It could be deduced from \eqref{DeltaphiH-3} that
	\begin{eqnarray}
		-\int_M d\phi(\nabla^{\nabla u}(\phi H))d\mu\geq \int_{M}\phi\Xi d\mu,
	\end{eqnarray}
	where $\Xi=\frac{1}{\phi}\{\mbox{The RHS of \eqref{DeltaphiH-3}}\}$. Suppose $\phi H$ achieves its positive maximum at $x_1$ in $B(p, 2R)$. By the same argument as in the compact case, one can prove that $\Xi\leq 0$ at $x_1$, otherwise $\phi H$ is a local subharmonic function in $B(p,2R)$, which provides a contradiction. Thus, considering that $\nabla^{\nabla u}(\phi H)=0$ at $x_1$, we know that
	\begin{eqnarray}\label{LapH-7}
		\begin{split}
			B\phi H\geq&\frac{2(1-\epsilon)}{N-\epsilon(N-n)}\frac{(\phi H)^2}{s^2q^2}+2\phi^2\frac{Ric^N(\nabla^{\nabla u}w)}{w^2}\\
			&+\phi^2\left[\frac{2(1-\epsilon)}{N-\epsilon(N-n)}\frac{(sq+s-1)^2}{s^2q^2}-2(\frac{1}{\epsilon}-1)\right]\frac{F_{\nabla u}^4(\nabla^{\nabla u}w)}{w^4}\\
			&+\phi\frac{4(1-\epsilon)(sq+s-1)}{N-\epsilon(N-n)s^2q^2}(\phi H)\frac{F_{\nabla u}^2(\nabla^{\nabla u}w)}{w^2}\\
			&+(2-6s)\phi^2w^{-\frac{2}{q}}\frac{F_{\nabla u}^2(\nabla^{\nabla u}w)}{w^2}-\frac2q(\phi H)g_{\nabla u}(\nabla^{\nabla u} \phi,\nabla^{\nabla u}\log w),
		\end{split}
	\end{eqnarray}
	is valid at $x_1$.
	Next we also consider \eqref{LapH-7} in two cases: (1) $C\leq 1$ and (2) $C>1$.\\
	
	Firstly, $u\leq 1$ when $C\leq 1$, so that $w\geq1$ and $1-w^{-\frac{2}{q}}\geq0$. Moreover, 
	\begin{eqnarray}\label{C1-RicciN}
		\frac{2Ric^N(\nabla^{\nabla u}w)}{w^2}\phi^2\geq -2K(2R)\phi^2\frac{F^2(\nabla^{\nabla u}w)}{w^2}\geq-2AK(2R)\phi H,
	\end{eqnarray}
	since $\phi\leq 1$ and $\frac{F_{\nabla u}^2(\nabla^{\nabla u}w)}{w^2}\leq H$. Furthermore, when $s=\frac{2}{3}$,
	\begin{eqnarray}\label{C1-est-1}
		&2w^{-\frac{2}{q}}\phi^2H+(2-6s)\phi^2\frac{F_{\nabla u}^2(\nabla^{\nabla u}w)}{w^2}w^{-\frac{2}{q}}
		\geq-(6s-4)\phi H=0.
	\end{eqnarray}

	Plugging \eqref{C1-RicciN} and \eqref{C1-est-1} into \eqref{LapH-7} and setting $0<q\leq \frac{3}{2(1+\sqrt{\frac{N-\epsilon(N-n)}{\epsilon}})}$ give that
	\begin{eqnarray}\label{ineq-B-1}
		\begin{split}
			B\phi H\geq &\frac{3(1-\epsilon)}{2q^2(N-\epsilon(N-n))}(\phi H)^2\\
			&-\frac{2}{q}(\phi H)g_{\nabla u}(\nabla^{\nabla u} \phi,\nabla^{\nabla u}\log w)-2AK(2R)\phi H.
		\end{split}
	\end{eqnarray}
	Adopting the Cauchy-Schwarz inequality and the estimates of the gradient of the cut-off function, we find that
	\begin{eqnarray}\label{C1-est-2}
		-\frac{2}{q}\phi Hg_{\nabla u}(\nabla^{\nabla u} \phi,\nabla^{\nabla u}\log w)\geq -\frac{2}{q}(\phi H)^{\frac{3}{2}}\frac{C_1}{R}.
	\end{eqnarray}
	Hence combining \eqref{ineq-B-1} and \eqref{C1-est-2} yields
	\begin{eqnarray}\label{ineq-B-2}
		B+2AK(2R)+\frac{2C_1}{qR}(\phi H)^{\frac12}\geq\frac{3(1-\epsilon)}{2q^2(N-\epsilon(N-n))}(\phi H).
	\end{eqnarray}
	Using the H\"older inequlity again gives that 
	$$\frac{2C_1}{qR}(\phi H)^{\frac12}\leq\frac{1-\epsilon}{2(N-\epsilon(N-n))q^2}(\phi H)+\frac{2(N-\epsilon(N-n))}{1-\epsilon}\frac{C_1^2}{R^2},$$ 
	so that
	\begin{eqnarray}
		\begin{split}
			\phi H&\leq \frac{q^2[N-\epsilon(N-n)]}{1-\epsilon}\left(B+\frac{2(N-\epsilon(N-n))}{1-\epsilon}\frac{C_1^2}{R^2}+2AK(2R)\right)\\
			&=\frac{q^2[N-\epsilon(N-n)]}{1-\epsilon}\left\{\frac{2C_1^2+C_2}{R^2}+2AK(2R)\right.\\
			&\quad+\frac{2(N+1-\epsilon(N+1-n))}{1-\epsilon}\frac{C_1^2}{R^2}\\
			&\quad\left.+\frac{C_1}{R}\left[\sqrt{K(2R)C(N,A)}\coth\left(R\sqrt{\frac{K(2R)}{C(N,A)}}\right)+C_0(K_0,A)\right]\right\}.
		\end{split}
	\end{eqnarray} 
	
	Secondly, when $C>1$, from the curvature condition and the H\"older inequality, we know that
	\begin{eqnarray}\label{C2-est-1}
		\begin{split}
			&\quad\phi^2\frac{2Ric^N(\nabla^{\nabla u}w)}{w^2}\geq-2AK(2R)\phi^2\frac{F_{\nabla u}^2(\nabla^{\nabla u}w)}{w^2}\\
			&\geq-\frac{2(1-\epsilon)(s-1)^2}{(N-\epsilon(N-n))s^2q^2}\phi^2\frac{F_{\nabla u}^4(\nabla^{\nabla u}w)}{w^4}-\frac{(N-\epsilon(N-n))s^2q^2}{2(1-\epsilon)(s-1)^2}\phi^2A^2K^2(2R),
		\end{split}
	\end{eqnarray}
	and 
	\begin{eqnarray}\label{C2-est-2}
		(6s-2)\phi^2\frac{F_{\nabla u}^2(\nabla^{\nabla u}w)}{w^2}w^{-\frac{2}{q}}\leq(3s-1)^2\phi^2\frac{F_{\nabla u}^4(\nabla^{\nabla u}w)}{w^4}+C^4\phi^2,
	\end{eqnarray}
	where we have used that $w^{-\frac{4}{q}}=u^4\leq C^4$. Meanwhile, the Cauchy-Schwarz inequality implies that
	\begin{eqnarray}\label{C2-est-3}
		\begin{split}
			-\frac{2}{q}\phi Hg_{\nabla u}(\nabla^{\nabla u}\phi,\nabla^{\nabla u}\log w)\leq&\frac{4(1-\epsilon)(sq+s-1)}{(N-\epsilon(N-n))s^2q^2}\frac{F_{\nabla u}^2(\nabla^{\nabla u}w)}{w^2}H\phi^2\\
			&+\frac{(N-\epsilon(N-n))s^2}{4(1-s)(sq+s-1)}\phi H\frac{F_{\nabla u}^2(\nabla^{\nabla u} \phi)}{\phi}.
		\end{split}
	\end{eqnarray}
	At last, we choose $s>1$ and $q>0$ (e.g. $s=2$ and $0\leq q\leq \frac{2\epsilon(1-\epsilon)}{(2+23\epsilon)(N-\epsilon(N-n))}$ for any chosen $\epsilon<1$) such that 
	$$\frac{2(1-\epsilon)}{(N-\epsilon(N-n))}\frac{s-1}{sq}\geq\frac{1}{\epsilon}-1+\frac{(3s-1)^2}{2}.$$
	We plug \eqref{C2-est-1}-\eqref{C2-est-3} into \eqref{LapH-7} and notice the above inequality to get
	\begin{eqnarray}\label{LapH-8}
		\begin{split}
			\frac{2(1-\epsilon)}{N-\epsilon(N-n)}\frac{(\phi H)^2}{s^2q^2}\leq&\left(\frac{(N-\epsilon(N-n))s^2C_1^2}{4(1-\epsilon)(sq+s-1)R^2}+B\right)(\phi H)\\
			&+\frac{(N-\epsilon(N-n))s^2q^2}{2(1-\epsilon)(s-1)^2}A^2K^2(2R)+C^4,
		\end{split}
	\end{eqnarray}
	where we have abandoned several nonnegative terms. Recall a result that if $a_0\leq a_1+a_0a_2$ for some $a_0,a_1,a_2\geq 0$, then $a_0\leq a_2+\sqrt{a_1}$. \eqref{LapH-8} implies that
	\begin{eqnarray}\label{LapH-9}
		\begin{split}
			\phi H\leq&\frac{(N-\epsilon(N-n))s^2q^2}{2(1-\epsilon)}\cdot\left\{\frac{(N-\epsilon(N-n))s^2C_1^2}{4(1-\epsilon)(sq+s-1)R^2}+\frac{(2C_1^2)}{R^2}+\frac{C_2}{R^2}\right.\\
			&\left.\quad+\frac{C_1}{R}\left[\sqrt{K(2R)C(N,A)}\coth\left(R\sqrt{\frac{K(2R)}{C(N,A)}}\right)+C_0(K_0,A)\right]\right\}\\
			&+\frac{(N-\epsilon(N-n))s^2q^2}{2(1-\epsilon)(s-1)}AK(2R)+sq\sqrt{\frac{N-\epsilon(N-n)}{2(1-\epsilon)}}C^2.
		\end{split}
	\end{eqnarray}
	Thus we finish the proof by taking $s=\frac23$ in \textit{Case 1.}
\end{proof}

\begin{rem}
	Theorem \ref{thm-noncpt-reduced} follows by taking $s=2$ in \textit{Case 2.} and $\epsilon=\frac12$, $q=\frac{2}{27(N+n)}$ in Theorem \ref{thm-noncpt-1}.
\end{rem}

At last, we give the proof of the Liouville theorem.

\begin{proof}[Proof of Theorem \ref{thm-Liouville}]
	In the case of $u\leq 1$, if $M$ is a compact Finsler manifold, Theorem \ref{thm-cpt-reduced} provides that
	\begin{eqnarray}\label{proof-lio}
		\frac{F^2(\nabla u)}{u^2}+\frac23(1-u^2)\leq0,
	\end{eqnarray}
	since $K=0$.  
	If $M$ is a forward complete Finsler manifold without boundary, Theorem \ref{thm-noncpt-reduced} also implies \eqref{proof-lio}, by letting $R$ tend to the infinity. However, \eqref{proof-lio} with $u\leq 1$ means that $\nabla u=0$ and $1-u^2=0$ in pointwise sense.
\end{proof}

\noindent

\bmhead{Acknowledgments}

The author is very grateful to the reviewers for their careful review and valuable comments. 

\section*{Declarations}

\begin{itemize}
	\item National Nature and Science foundation of China (No.12001099, 12271093)
	\item Authors have no conflict of interest.
	\item Ethics approval is not applicable
	\item Data sharing is not applicable
\end{itemize}

\end{document}